# ROTATION PREVENTS FINITE-TIME BREAKDOWN

HAILIANG LIU AND EITAN TADMOR

## Contents

ABSTRACT. We consider a two-dimensional convection model augmented with the rotational Coriolis forcing, $U_t + U \cdot \nabla_x U = 2kU^\perp$, with a fixed $2k$ being the inverse Rossby number. We ask whether the action of dispersive rotational forcing alone, $U^\perp$, prevents the generic finite time breakdown of the free nonlinear convection. The answer provided in this work is a conditional yes. Namely, we show that the rotating Euler equations admit global smooth solutions for a subset of generic initial configurations. With other configurations, however, finite time breakdown of solutions may and actually does occur. Thus, global regularity depends on whether the initial configuration crosses an intrinsic, $\mathcal{O}(1)$ critical threshold, which is quantified in terms of the initial vorticity, $\omega_0 = \nabla \times U_0$, and the initial spectral gap associated with the $2 \times 2$ initial velocity gradient, $\eta_0 := \lambda_2(0) - \lambda_1(0)$, $\lambda_j(0) = \lambda_j(\nabla U_0)$. Specifically, global regularity of the rotational Euler equation is ensured if and only if $4k\omega_0(\alpha) + \eta_0^2(\alpha) < 4k^2$, $\forall \alpha \in \mathbb{R}^2$.

We also prove that the velocity field remains smooth if and only if it is periodic. An equivalent Lagrangian formulation reconfirms the critical threshold and shows a global periodicity of velocity field as well as the associated particle orbits. Moreover, we observe yet another remarkable periodic behavior exhibited by the *gradient* of the velocity field. The spectral dynamics of the Eulerian formulation, [18], reveals that the vorticity and the eigenvalues (and hence the divergence) of the flow evolve with their own path-dependent period. We conclude with a kinetic formulation of the rotating Euler equation.

## 1. Introduction and statement of main results

Finite-time breakdown is a familiar trademark of nonlinear convection mechanism. Consider the canonical example of an $N$-dimensional system of free transport equations, $U_t + U \cdot \nabla_x U = 0$. It follows — consult Corollary 2.2 below, that the solution $U(t, \cdot)$ will lose its initial regularity at a finite-time if and only if an eigenvalue of the initial velocity gradient crosses the negative real axis, i.e., iff there exists at least one eigenvalue, $\lambda(0, x) := \lambda(\frac{\partial U_i(0,x)}{\partial x_j})$, such that $\lambda(0, x) \in \mathsf{IR}^-$. Consequently, finite-time breakdown is a generic phenomenon for the free nonlinear transport. Thus, for example, irrotational initial data $\nabla_x \times U(0, x) = 0$ — where all eigenvalues $\lambda_j(t, x)$ remain real, will necessarily lead to finite-time breakdown, except for non-generic cases where $\lambda_j(0, x) \geq 0$, $\forall j, x$, requiring, in particular, that the initial divergence is *globally* positive, $\nabla_x \cdot U(0, x) > 0$. This general $N$-dimensional scenario is completely analogous to the 1D inviscid Burgers' equation, $U_t + UU_x = 0$, where solutions of the latter will necessarily reach a finite-time breakdown except for the non-generic case of monotonically increasing initial data.

Physically relevant models are governed by the fundamental Eulerian convection equation augmented by proper forcing $F$,

$$(1.1) \qquad U_t + U \cdot \nabla_x U = F.$$

Here, there is a competition between the finite-time breakdown dynamics driven by nonlinear convection and the balancing act of nonlinear forcing, $F$. Different models show up in different contexts dictated by the different modeling of such forcing. Three prototypes are dissipation, relaxation and dispersion. It is well known that if (1.1) is augmented with a sufficiently large amount of either dissipation or relaxation, then (1.1) admits a global smooth solution for a rich enough class of initial data. In both cases of dissipation and relaxation, global existence is secured by enforcing a sufficiently large amount of energy decay. Dispersive forcing, however, is different. The dispersive KdV equation, for example, $U_t + UU_x = U_{xxx}$ is a case in point. It admits global smooth solution while keeping the $L^2$-energy invariant in time. In this paper, we study the regularity of the 2D convection model augmented by rotational forcing,

$$(1.2) \qquad U_t + U \cdot \nabla_x U = 2kJU, \qquad J := \begin{pmatrix} 0 & 1 \\ -1 & 0 \end{pmatrix},$$

subject to initial conditions, $U(0, x) = U_0(x)$. Here $2k = \epsilon^{-1}$ where $\epsilon$ is the Rossby number [22], $\epsilon = \frac{\overline{U}}{2\Omega \overline{L}}$, expressed in terms of the characteristic length $\overline{L}$, characteristic speed, $\overline{U}$, and the amplitude of angular velocity $\Omega$ of the rotating body, consult [14, 21]. With these parameters the system evolves on a characteristic time scale $t \sim \overline{L}/\overline{U}$.

The system admits a global energy invariant in time, which is independent of the amplitude of rotation encoded by the constant $k$ on the RHS of (1.2). To see this, we note that (1.2) is formally equivalent to the extended $3 \times 3$ system,

$$(1.3) \qquad \partial_t \rho + \nabla_x \cdot (\rho U) = 0, \quad x \in \mathsf{IR}^2, \quad t \in \mathsf{IR}^+,$$
$$(1.4) \qquad \partial_t(\rho U) + \nabla_x \cdot (\rho U \otimes U) = 2k\rho JU,$$

which are the usual statements of conservation of mass and Newton's second law, governing the local density $\rho = \rho(t, x)$ and the velocity field $U := (u, v)(t, x)$. The usual manipulation, $-1/2|U|^2 \times (1.3) + U^\top \times (1.4)$ and the skew-symmetry form induced by the rotational forcing imply

$$\partial_t \left(\frac{1}{2}\rho|U|^2\right) + \nabla_x \cdot \left(\frac{1}{2}\rho U|U|^2\right) = 2k\rho \langle U, JU \rangle = 0.$$

The global invariance of the energy follows

$$E(t) := \frac{1}{2} \int_x \rho(t,x)|U(t,x)|^2 dx = E(0).$$

The system (1.3)-(1.4) coincides with a simplified version of the 2D shallow-water equations, where the additional pressure terms are ignored. The only remaining forcing is the rotational Coriolis forcing, and our main quest in this paper is whether the action of dispersive rotational forcing alone prevents the generic finite time breakdown of nonlinear convection. The answer, outlined in Section 4 is a conditional yes. Namely, we show that (1.2) admits global smooth solutions for a subset of generic initial configurations, $U_0$. With other initial configurations, however, the finite time breakdown of solutions may – and actually does occur. Thus, global regularity depends on whether the initial configuration crosses an intrinsic, $\mathcal{O}(1)$ critical threshold, which is quantified in terms of the initial vorticity, $\omega_0 := \nabla_x \times U_0$ and the initial spectral gap, $\Gamma_0 := (\lambda_2(0) - \lambda_1(0))^2$.

**Theorem 1.1** (Critical threshold for rotation forcing). *Consider the 2D rotational flow (1.2) with $k > 0$. Then the solution of (1.2) with initial data $U(x,0)\big|_{x=\alpha} = U_0(\alpha)$ remains smooth for all time, $-\infty < t < \infty$, if and only if the initial data $U_0$ satisfy*

(1.5) $$i_0(\alpha) := 4k[k - \omega_0(\alpha)] - \Gamma_0(\alpha) > 0, \quad \forall \alpha \in \mathbb{R}^2.$$

*Moreover, as the smooth solution evolves along its particle path, then the vorticity, $\omega(t) = \omega(t,\alpha)$ and the eigenvalues $\lambda_j$ (and hence the divergence, $d(t) = div_x U(t,\alpha)$) form a periodic orbit in phase space, with a path-dependent period, $\overline{T} = \overline{T}(\alpha)$, given by*

(1.6) $$\overline{T} = \frac{2}{k} \int_{-\pi/2}^{\pi/2} \frac{d\theta}{(\theta_0^{-1} + \theta_0) + (\theta_0^{-1} - \theta_0)\sin\theta}.$$

*Here $\theta_0 = \theta(\alpha) < 1$ is determined by the initial data*

(1.7) $$\theta_0 = \frac{\sqrt{1 + 8kp_0} - 1}{\sqrt{1 + 8kp_0} + 1}, \quad p_0 := \frac{\sqrt{i_0}}{d_0^2 + (\sqrt{i_0} - 2k)^2}.$$

Several remarks are in order.

1. We note that the critical threshold phenomena is independent of the initial divergence $d_0 := \text{div}_x U_0$.

2. Let us point out that system (1.2) could be viewed as a crossroad between the 2D shallow-water equations and the so-called 2D pressureless equations, e.g. [1, 2, 3, 4, 8, 13], corresponding to Theorem 1.1 with $k = 0$,

(1.8) $$\partial_t \rho + \nabla_x \cdot (\rho U) = 0, \quad x \in \mathbb{R}^2, \quad t \in \mathbb{R}^+,$$

(1.9) $$\partial_t(\rho U) + \nabla_x \cdot (\rho U \otimes U) = 0.$$

According to Theorem 1.1, the pressureless system admits a global smooth solution forward (respectively, reversible) in time, if and only if $\lambda_j(0) \notin \mathbb{R}^-$, respectively iff $\lambda_j(0) \notin \mathbb{R}$. The latter is equivalently expressed in Theorem 1.1 requiring $\Gamma_0 := (\lambda_2(0) - \lambda_1(0))^2 < 0$.

3. In particular, (1.2) does admit global smooth solutions with negative initial divergence in contrast to the free transport ($k = 0$) equation discussed in the introduction. It follows that rotation prevents finite time breakdown, either by a large Coriolis forcing ($k >> 1$) or a large initial rotation ($\Gamma_0 << 0$).

4. If we set $y$-independent initial data, then (1.2) is reduced to the one-dimensional system

$$\begin{aligned} u_t + uu_x &= 2kv, \\ v_t + uv_x &= -2ku, \end{aligned}$$

with critical threshold $(u_0'(\alpha))^2 - 4kv_0'(\alpha) < 4k^2$. To interpret Theorem 1.1 in this simplified setting, we observe that the gradient $(\omega, d) := (-v_x, u_x)$ solves a coupled system

$$(1.10) \qquad (\partial_t + u\partial_x)\omega + d\omega = 2kd,$$

$$(1.11) \qquad (\partial_t + u\partial_x)d + d^2 = -2k\omega,$$

and a straightforward computation reveals the path-dependent invariant along the particle path, $\dot{X}(t) = u(X, t), X(0) = \alpha$,

$$\frac{(2k - \omega)^2}{d^2 + \omega^2} = B_0, \quad B_0 = B_0(\alpha) := \frac{(2k - \omega_0(\alpha))^2}{d_0^2(\alpha) + \omega_0(\alpha)^2}.$$

The critical threshold statement in this case reads $B_0 > 1$, stating that the gradient $(\omega, d)$ forms a closed elliptical orbit in the phase plane (whereas $B_0 \leq 1$ corresponds to unbounded parabolic/hyperbolic orbits). Following the analysis in Section 4, we also obtain a path-dependent period for the gradient

$$\overline{T} = \frac{2}{k} \int_{-\pi/2}^{\pi/2} \frac{d\theta}{(\theta_0^{-1} + \theta_0) + (\theta_0^{-1} - \theta_0)\sin\theta}, \quad \theta_0 = \frac{\sqrt{B_0} - 1}{\sqrt{B_0} + 1}.$$

Such path-dependent period of the gradient reflects the fact that its governing system (1.10),(1.11), is a nonlinear perturbation of the harmonic oscillator. As the Rossby number approaches zero, however, $k >> 1$, $\theta_0 \sim 1$, and the above path-dependent period $\overline{T}$ is approaching the global inertial period $\pi/k$ (the harmonic oscillator period).

5. As we shall see in Theorem 1.2 below, sub-critical initial data yield smooth velocity fields, $U(X(t))$, with time period $\overline{T} = \pi/k$, or — expressed in terms of the original non-scaled time units, a period $T = \overline{T}\,\overline{U}/\overline{L} = \pi/\Omega$. This period of the particle orbits is related to the global revolution of the plane. Theorem 1.1 points out yet another remarkable property for a portion of the gradient of $U(t, \cdot)$, namely, the divergence $d(t, \alpha)$, and the vorticity $\omega(t, \alpha)$ which exhibit a local, path-dependent period dictated by the unique initial parameter, $8kp_0$. It is instructive to compute the period predicted in Theorem 1.1, using configurations similar to those encountered in various applications. Let us illustrate a couple of examples taken from [14]. For the Gulf Stream, with the Rossby number $\epsilon = 0.07$, $\overline{L} = 100$ km and $\overline{U} = 1m/sec$, we find that the vorticity and divergence of the flow keep repeating themselves every $T = \overline{T}\,\overline{L}/\overline{U} \sim 11.7hrs$; for the weather system we have $\epsilon = 0.14$ with $\overline{L} = 1000km$, $\overline{U} = 20m/sec$, and the vorticity/divergence exhibit a period of $T = \overline{T}\,\overline{L}/\overline{U} \sim 12.2hrs$. It is also interesting to see how this path-dependent gradient period be influenced by the small Rossby numbers. After rescaling we may assume initial configuration such that $d_0 \sim \omega_0 \sim 1$, for which a small Rossby number yields $i_0 \sim 4k^2$, $p_0 \sim 2k$ and hence $\theta_0 \sim 1$. Restored in terms of the original time scale, $T = \overline{T}\,\overline{L}/\overline{U}$, the period is given by the $\theta_0$-dependent elliptic integral,

$$T = \frac{2}{\Omega} \int_{-\pi/2}^{\pi/2} \frac{d\theta}{(\theta_0^{-1} + \theta_0) + (\theta_0^{-1} - \theta_0)\sin\theta} \sim \frac{\pi}{\Omega},$$

which is close to the inertial period when the Rossby number is small. For the earth core, for example, we have $T = 11.95hrs$ with $\epsilon = 2 \times 10^{-7}$, ($\overline{L} = 3000km$ and $\overline{U} = 0.1cm/sec$), whereas for Jupiter's Red Spot we have a Rossby number $\epsilon = 0.015$ (with $\overline{L} = 10^4 km$ and $\overline{U} = 0.1cm/sec$) and the velocity gradient period $T \sim 5.13hrs$. We should point out the difference between the period of the velocity field vs. the velocity gradient periods, which is due to our tracking of the flow dynamics along the particle path.

Of course, one should not expect the current cartoon model to provide a faithful description of the full model, as other forces ignored at this stage — magnetic forces, pressure, etc, will play the important role. Nevertheless, the above rough approximations are interesting for their own sake, in particular, since the flow is predicted to be periodic once smooth solutions are secured for subcritical initial data. Surprisingly, the periods computed above fall within the physical range. It will be challenging to refine the estimate, by taking into account the other forces which augment with the rotation model.

To put our study in a proper perspective we recall that there has been a considerable amount of literature available on the global behavior of nonlinear convection driven by rotational forcing and related problems, from rotating shallow-water model [12, 15, 23] to rotating incompressible Euler and Navier-Stokes equations [5, 6, 9, 7]. The common feature of the flows studied in this context are rotation dominated flow for which the Rossby number $\epsilon$ is small. It is well known that large-scale atmospheric (or oceanic) fields are in permanent process of Rossby (or geostrophic) adjustment [20]. The flow structure has been extensively studied in terms of $\epsilon$, say in [5, 12], based on the averaging of the interaction of the fast waves of the rotating Euler equation, two dimensional structures were shown to emerge in the limit $\epsilon \to 0$; for bounds for the vertical gradients of the Lagrangian displacement that vanish linearly with the maximal local Rossby number [9]; for a nonlinear theory of geostrophic adjustment for the rotating shallow-water model under the assumption of the smallness of the Rossby number [23]; consult [15] for the analysis of an approximation of the rotating shallow-water equation.

When dealing with the questions of time regularity for Eulerian dynamics without damping, one encounters several limitations with the classical stability analysis. Among other issues, we mention that

(i) the stability analysis does not tell us how large perturbations are allowed before losing stability – indeed, the smallness of the initial perturbation is essential to make the energy method work, e.g. the 3D Navier-Stokes equation [16];

(ii) the steady solution may be only conditionally stable due to the weak dissipation in the system, say in the 1D Euler-Poisson equations [11].

To address these difficulties we advocated, in [11, 17, 18], a new notion of critical threshold (CT) which describes conditional stability, where the answer to the question of global vs local existence depends on whether the initial configuration crosses an intrinsic, $\mathcal{O}(1)$ critical threshold. Little or no attention has been paid to this remarkable phenomena, and our goal is to bridge the gap of previous studies on the behavior of rotational Euler equations, a gap between the regularity of Eulerian solutions in the small and their finite-time breakdown in the large. The critical threshold (CT) was completely characterized for the 1D Euler-Poisson system in terms of the relative size of the initial velocity slope and the initial density; consult [17, 24] for the CT for the convolution model for conservation laws; Moving to the multi-D setup, one has first to identify the proper quantities which govern the critical threshold phenomena. In [19] we have shown that these quantities depend in an essential manner on the *eigenvalues* of the velocity gradient matrix, $\lambda(\nabla_x U)$.

The critical threshold for the current rotation model can be also obtained, in a straightforward manner, through a Lagrangian flow formulation. This is summarized in the Theorem 1.2 below. We should point out that it was the spectral dynamics analysis of $\lambda(\nabla_x U)$ that led us to the CT formulation in the first place, which in turn was then sought within Lagrangian formulation. In Section 5 we prove

**Theorem 1.2** (Flow map for rotation forcing). *The flow map,* $\dot{X}_\alpha := \frac{dX_\alpha}{dt} = U(X_\alpha)$, $X_\alpha(0) = \alpha$ *associated with (1.2),* $\ddot{X}_\alpha = 2kJ\dot{X}_\alpha$ *is given by*

$$X_\alpha(t) = \frac{1}{2k}J^{-1}e^{2kJt}U_0(\alpha) + \alpha - \frac{1}{2k}J^{-1}U_0(\alpha).$$

*For sub-critical initial data, (1.5), this flow map is invertible and periodic with an inertial period* $\overline{T} = \pi/k$. *The velocity field* $U(X_\alpha(t)) = U_0(\alpha) + 2kJ(X_\alpha(t) - \alpha)$ *shares the same inertial period.*

At this point, one may wonder whether this inertial period is none other than the planet rotation. Actually the two are not the same; the rotating plane completes one revolution in a time equal to $\frac{2\pi}{\Omega}$, while the period of the particle path expressed in the original non-scaled variables is $T = \overline{T}L/\overline{U} = \frac{\pi}{\Omega}$. Thus, the particle goes around its orbit twice as the plane accomplishes a single revolution, which is consistent with the observation in [10].

Finally we conclude in Section 6 with a kinetic formulation of the current rotation model.

**Theorem 1.3.** *The rotation model (1.2) admits for the following kinetic formulation*

$$\partial_t f + \xi \cdot \nabla_x f + 2kJ\xi \cdot \nabla_\xi f = \frac{1}{\epsilon}(M - f),$$

*where* $M_{\{\rho, U\}}(\xi)$ *is the Maxwellian given by*

$$M = \frac{\rho}{\sqrt{\pi T}}e^{-|\xi - U|^2/T}, \quad \xi = (\xi_1, \xi_2) \in \mathrm{IR}^2,$$

*where $\rho$ and $U$ are macroscopic density and velocity respectively, and $T$ being an arbitrary fixed temperature.*

In Section 4 and 5 below, we quantify the same critical threshold using the Eulerian and Lagrangian formulations, and it would be of interest to derive the same critical threshold directly using the kinetic formulation in Theorem 1.3.

## 2. Spectral dynamics

We consider a general nonlinear transport equation (1.1), $U_t + U \cdot \nabla_x U = F$, and we trace the evolution of $\nabla_x U$ in terms of its eigenvalues, $\lambda := \lambda(\nabla_x U)(t, x)$. The following result is in the heart of matter.

**Lemma 2.1.** [ Spectral dynamics, [18, Lemma 3.1] ] . *Let* $\lambda := \lambda(\nabla_x U)(t, x)$ *denote an eigenvalue of* $\nabla_x U$ *with corresponding left and right normalized eigenpair,* $\langle \ell, r \rangle = 1$. *Then $\lambda$ is governed by the forced Riccati equation*

$$\partial_t \lambda + U \cdot \nabla_x \lambda + \lambda^2 = \langle \ell, \nabla_x F r \rangle.$$

As an immediate corollary we obtain the precise description for finite time breakdown of free nonlinear transport.

**Corollary 2.2.** [Finite time breakdown of free transport, [18, Lemma 4.1]]. *The free nonlinear $N$-dimensional transport*

$$\partial_t U + U \cdot \nabla_x U = 0, \quad x \in \mathbb{R}^N, \tag{2.1}$$

*admits global smooth solution forward in time, $t > 0$, if and only if the eigenvalues of its initial velocity gradient, $\lambda := \lambda(\nabla_x U)$, satisfy $\lambda(0, x) \notin \mathbb{R}^-$. Likewise, it admits a globally smooth, time-reversible solution for $-\infty < t < \infty$ iff $\lambda(0, x) \notin \mathbb{R}$.*

For the proof, we note that the eigenvalues, governed by the homogeneous Riccati equations, propagate along the particle path $x = x(t, \alpha)$,

$$\lambda(t, x) = \frac{\lambda(0, \alpha)}{t\lambda(0, \alpha) + 1}.$$

We note in passing that the rotational system (1.2) is to the full shallow-water equations as the free transport model (2.1) is to the full Euler equations. The existence of a critical threshold phenomena associated with global linear forcing model was first identified by us [18], although the exact configuration cannot be obtained in such generality. The current paper provides a precise description of the critical threshold for the 2D rotational system (1.2). In particular, we use the Spectral Dynamics Lemma to obtain remarkable explicit formulae for the critical threshold surface summarized in the main Theorems 1.1. Taking the gradient of the velocity equation (1.2), we find that the velocity gradient field, $\nabla_x U$, solves the following matrix equation

$$\partial_t(\nabla_x U) + U \cdot \nabla_x(\nabla_x U) + (\nabla_x U)^2 = 2kJ\nabla_x U. \tag{2.2}$$

Using the spectral dynamics Lemma 2.1, we obtain the spectral dynamics equations

$$\partial_t \lambda_1 + U \cdot \nabla_x \lambda_1 + \lambda_1^2 = 2k\lambda_1 \langle l_1, Jr_1 \rangle, \tag{2.3}$$

$$\partial_t \lambda_2 + U \cdot \nabla_x \lambda_2 + \lambda_2^2 = 2k\lambda_2 \langle l_2, Jr_2 \rangle, \tag{2.4}$$

where $\lambda_i, i = 1, 2$ are eigenvalues of the velocity gradient field $\nabla_x U$ associated with left (row) eigenvectors $l_i$ and right (column) eigenvectors $r_i$. Since $J$ is skew-symmetric we have $Jr_1 = \alpha_1 l_2^\top$ and $Jr_2 = \alpha_2 l_1^\top$. Noting that $l_2 r_2 = l_1 r_1 = 1$, one then has $\alpha_1 = \langle r_2, Jr_1 \rangle$ and

$$\alpha_2 = \langle r_1, Jr_2 \rangle = -\langle r_2, Jr_1 \rangle = -\alpha_1.$$

Therefore, (2.3)-(2.4) now read

$$\partial_t \lambda_1 + U \cdot \nabla_x \lambda_1 + \lambda_1^2 = 2k\lambda_1 \langle r_2, Jr_1 \rangle \langle l_1, l_2 \rangle, \tag{2.5}$$

$$\partial_t \lambda_2 + U \cdot \nabla_x \lambda_2 + \lambda_2^2 = -2k\lambda_2 \langle r_2, Jr_1 \rangle \langle l_1, l_2 \rangle, \tag{2.6}$$

from which we deduce that the spectral gap $\eta := \lambda_2 - \lambda_1$ and divergence $d := \lambda_2 + \lambda_1$, satisfy

$$\partial_t \eta + U \cdot \nabla_x \eta + d\eta = -2kd\langle r_2, Jr_1 \rangle \langle l_1, l_2 \rangle$$

and

$$\partial_t d + U \cdot \nabla_x d + \frac{d^2 + \eta^2}{2} = -2k\eta \langle r_2, Jr_1 \rangle \langle l_1, l_2 \rangle. \tag{2.7}$$

On the other hand, differentiation of (1.2) yields the $\nabla_x U$-equation (2.2), i.e.,

$$(2.8) \quad (\partial_t + U \cdot \nabla_x)\begin{pmatrix} u_x & u_y \\ v_x & v_y \end{pmatrix} + \begin{pmatrix} u_x^2 + v_x u_y & du_y \\ dv_x & v_x u_y + v_y^2 \end{pmatrix} = 2k\begin{pmatrix} v_x & v_y \\ -u_x & -u_y \end{pmatrix},$$

which in turn – using the LHS of (2.5)-(2.6) to express $\lambda_1^2 + \lambda_2^2 \equiv (d^2 + \eta^2)/2$, leads to

$$(2.9) \quad \partial_t d + U \cdot \nabla_x d + \frac{d^2 + \eta^2}{2} = -2k\omega,$$

$$(2.10) \quad \partial_t \omega + U \cdot \nabla_x \omega + d\omega = 2kd.$$

Equating the expressions on the right of (2.9) and (2.7) we find

$$-2k\omega = -2k\eta \langle r_2, Jr_1\rangle \langle l_1, l_2\rangle.$$

Thus, the scaled product of the eigenvectors measures the ratio of vorticity over the spectral gap in following manner

$$(2.11) \quad \langle r_2, Jr_1\rangle \langle l_1, l_2\rangle = \frac{\omega}{\eta}.$$

When the spectral gap $\eta$ shrinks to zero, the scaled product becomes unbounded due to the degeneracy of eigenvectors. When the vorticity $\omega$ shrinks to zero, (2.11) recovers the symmetry of $\nabla_x U$ which is reflected through the orthogonality of $\ell_1$ and $\ell_2$. Equipped with the above relations we come up with a closed system for $(\omega, d, \eta)$ along the particle path (here and below $' \equiv \partial_t + U \cdot \nabla_x$)

$$\omega' + d\omega = 2kd,$$

$$d' + \frac{d^2 + \eta^2}{2} = -2k\omega,$$

$$\eta' + d\eta = -2k\frac{d\omega}{\eta}.$$

Note that the spectral gap may become purely imaginary when eigenvalues are complex. To avoid the discussion on the complex solution of the above system, we introduce the following real variable

$$\Gamma := \eta^2.$$

Using the above equations we have

$$\Gamma' = 2\eta\eta' = 2d[-2k\omega - \Gamma].$$

Note that the sign of $2k$ indicates the direction of the rotational forcing, and the vorticity measures the rotation in the flow. In order to combine these two effects we introduce $\varphi := 4k^2 - 2k\omega$, and thus obtain a closed system for $W := (\varphi, d, \Gamma)^\top$

$$(2.12) \quad \varphi' = -d\varphi,$$

$$(2.13) \quad d' = -\frac{d^2 + \Gamma}{2} + \varphi - 4k^2,$$

$$(2.14) \quad \Gamma' = 2d[\varphi - 4k^2 - \Gamma].$$

We shall use this system to describe the dynamics of the velocity gradient field. Linearization of the above system around $W^* = (\varphi^*, d^*, \Gamma^*)^\top$ gives the linear system $W' = A(W^*)(W - W^*)$ with

$$A = \begin{pmatrix} -d^* & -\varphi^* & 0 \\ 1 & -d^* & -\frac{1}{2} \\ 2d^* & 2\varphi^* & -2d^* \end{pmatrix}$$

The corresponding eigenvalues of $A$ at critical points $(\varphi^*, 0, \Gamma^*)$ are $\lambda_1 = 0$, $\lambda_{2,3} = \pm\sqrt{-2\varphi^*}$. The classical stability analysis based on linearization is not sufficient to predict the global time dynamics.

## 3. Material Invariants

It follows from equations (2.12) and (2.14) we obtain

$$\frac{d\varphi}{d\Gamma} = \frac{-\varphi}{2(\varphi - 4k^2 - \Gamma)},$$

which upon integration gives the first material invariant

$$(3.1) \qquad \frac{2\varphi - \Gamma - 4k^2}{\varphi^2} = \mathcal{C}_0, \qquad \mathcal{C}_0 \equiv \mathcal{C}_0(\alpha) := \frac{2\varphi_0 - \Gamma_0 - 4k^2}{\varphi_0^2}.$$

This material invariant enables us to reduce the full system (2.12)-(2.14) to the following system

$$(3.2) \qquad \varphi' = -\varphi d,$$

$$(3.3) \qquad d' = -\frac{1}{2}[d^2 + 4k^2 - \mathcal{C}_0\varphi^2].$$

In order to have global bounded solution it is necessary to assume $\mathcal{C}_0(\alpha) > 0$, i.e.,

$$(3.4) \qquad \Gamma_0 < 2\varphi_0 - 4k^2 \equiv 4k(k - \omega_0),$$

for otherwise, (3.3) will be majored by the Riccati equation $d' \leq -\frac{d^2}{2} - 2k^2$, which would lead to the finite time breakdown. As we shall see below in §4, the positivity of $C_0(\alpha)$ is also sufficient for the global bounded solution.

Another material invariant is obtained in the following manner. Following [19, Lemma 2.2], we set $q = d^2$ which yields

$$\frac{dq}{d\varphi} = 2d\frac{d'}{\varphi'} = \frac{q + 4k^2 - \mathcal{C}_0\varphi^2}{\varphi}.$$

Integration yields a second material invariant

$$(3.5) \qquad \frac{d^2 + 4k^2 + \mathcal{C}_0\varphi^2}{\varphi} = \mathcal{D}_0, \qquad \mathcal{D}_0 \equiv \mathcal{D}_0(\alpha) := \frac{d_0^2 + 4k^2 + \mathcal{C}_0\varphi_0^2}{\varphi_0};$$

and together with (3.1) we end up with a second independent material invariant

$$\frac{d^2 - \Gamma + 2\varphi}{\varphi} = \mathcal{D}_0(\alpha).$$

In summary we have

**Lemma 3.1.** *Let $\varphi := 4k^2 - 2k\omega$, $d := divU = tr(\nabla_x U)$ and $\Gamma := (\lambda_2 - \lambda_1)^2$ be the solution of the dynamical system (2.12)-(2.14), associated with the rotational system (1.2). Then we have the following material invariants along particle path $(t, X_\alpha(t))$,*

$$\left.\frac{2\varphi - \Gamma - 4k^2}{\varphi^2}\right|_{(t,X_\alpha(t))} = \frac{2\varphi_0(\alpha) - \Gamma_0(\alpha) - 4k^2}{\varphi_0^2(\alpha)}, \tag{3.6}$$

$$\left.\frac{d^2 - \Gamma}{\varphi}\right|_{(t,X_\alpha(t))} = \frac{d_0^2(\alpha) - \Gamma_0(\alpha)}{\varphi_0(\alpha)}. \tag{3.7}$$

## 4. Critical Thresholds

As we observed earlier, the positivity of condition, $\mathcal{C}_0(\alpha) > 0$, is necessary for global bounded solution, for otherwise

$$d' < -\frac{1}{2}[4k^2 + d^2],$$

which would imply that $d$, and hence $\varphi$, become unbounded in a finite time. We shall show that the same positivity condition, $\mathcal{C}_0(\alpha) > 0$, is in fact also sufficient for the existence of global bounded solution. For $\mathcal{C}_0 > 0$, the reduced system (3.2)-(3.3) has two unique equilibrium points in the finite plane, $(\varphi_\pm^*, d) = (\pm\frac{2k}{\sqrt{\mathcal{C}_0}}, 0)$. The local behavior of the solution depends on the properties of these critical points. We note that since $\varphi = 0$ is an invariant set, then $\varphi_0\varphi(t) > 0$ for all time, and we therefore concentrate on the solution behavior for $\varphi_0 > 0$, with the other case of $\varphi_0 < 0$ being handled similarly. On the right plane $\varphi > 0$, the coefficient matrix of linearized system of (3.2)-(3.3) around the equilibrium point $(\varphi_+^*, d) = (\frac{2k}{\sqrt{\mathcal{C}_0}}, 0)$ is

$$\begin{pmatrix} 0 & -\varphi_+^* \\ \mathcal{C}_0\varphi_+^* & 0 \end{pmatrix},$$

with purely imaginary eigenvalues, $\pm(\sqrt{\mathcal{C}_0}\varphi_+^*)i$. This means that the bounded trajectory is possibly a periodic solution or limit circle. Observe that if $(\varphi(t), d(t))$ is a solution to (3.2), (3.3), so is $(\varphi(-t), -d(-t))$. Such symmetry implies that $(\varphi_+^*, 0)$ is a center and the trajectory in a neighborhood of this equilibrium point is periodic.

In order to clarify the global behavior of the flow around the center, we appeal to the material invariant (3.5) which we rewrite as

$$V_+(\varphi, d) := \frac{d^2 + (\sqrt{\mathcal{C}_0}\varphi - 2k)^2}{\varphi}.$$

$V_+(\cdot)$ is a positive definite Liapunov function for $\varphi > 0$, and achieves its global minimum, $V_+ = 0$, at the equilibrium point $(\varphi_+^* = \frac{2k}{\sqrt{\mathcal{C}_0}}, 0)$. A family of closed orbits in the phase plane $(\varphi, d)$ can be expressed as the level set curve

$$V_+(\varphi, d) = Const > 0,$$

since $V_+$ is material invariant in the sense that $\frac{dV_+}{dt} = 0$. Similarly, on the left plane $\varphi < 0$, one may use the Liapunov functional

$$V_-(\varphi, d) = \frac{d^2 + (\sqrt{\mathcal{C}_0}\varphi + 2k)^2}{-\varphi},$$

whose level set curve determines a family of closed orbit on the left plane centered around $(\varphi_-^* = -\frac{2k}{\sqrt{C_0}}, 0)$. The global behavior of the solutions is summarized in

**Lemma 4.1.** [Bounded solutions are periodic]
 $\{i\}$ Given initial data $(\varphi_0, d_0, \Gamma_0)$ for the system (2.12)-(2.14), it admits a global bounded solution if and only if the initial data lie in the sub-critical region where

(4.1) $$\Gamma_0 < 2\varphi_0 - 4k^2.$$

$\{ii\}$ Bounded solutions of (2.12)-(2.14) are necessarily periodic. The periodic orbit on the right plane $\varphi > 0$ lies on the ellipse $d^2 + (\sqrt{C_0}\varphi - 2k)^2 = (\mathcal{D}_0 - 4k\sqrt{C_0})\varphi$, with $C_0$ and $\mathcal{D}_0$ are determined by the initial data, (3.1) and (3.5).

*Proof.* By (4.1), the initial data for the bounded orbit must satisfy $\Gamma_0 < 2\varphi_0 - 4k^2$, i.e. $C_0 > 0$. From the fact that $V(\varphi, d) = Const > 0$ is the material invariant, we find that $\Gamma$ can be expressed as
$$d^2 + (\sqrt{C_0}\varphi - 2k)^2 = (\mathcal{D}_0 - 4k\sqrt{C_0})\varphi.$$
This is an ellipse, along which the corresponding solution $(\varphi = 4k^2 - 4k\omega, d)$ being periodic. Moreover, the invariants (3.6), (3.7) imply that $\Gamma$ shares the same period along the path. It follows that the boundedness of the divergence implies the boundedness of the whole velocity gradient [19][Lemma 2.1]: the anti-trace of (2.8) $r := v_x + u_y$ satisfies $r' + rd = -2ks$ followed by a $s := u_x - v_y$-equation, $s' + sd = 2kr$ yield the boundedness of $\nabla U$. In fact solving the above two equations in terms of the divergence $d$ we obtain

$$\frac{r}{s} = \tan(\tan^{-1}(r_0/s_0) - 2kt), \quad r^2 + s^2 = (r_0^2 + s_0^2)\exp\left(-2\int_0^t d(\xi)d\xi\right).$$

By the periodicity of $d$ and its symmetry about the axis $d = 0$ we see that $s^2 + r^2$ shares the same period with $d$; This is also clear from the identity $s^2 + r^2 = \Gamma + \omega^2$. A combination of the above facts lead to

$$r = \sin(\tan^{-1}(r_0/s_0) - 2kt)\sqrt{\Gamma + \omega^2}, \quad s = \cos(\tan^{-1}(r_0/s_0) - 2kt)\sqrt{\Gamma + \omega^2}.$$

As a product of two periodic functions with period $\pi/k$ and $\overline{T}$ respectively, the anti-trace $r$, also the anti-vorticity $s$, is periodic only if the ratio of these two periods, i.e.,

$$2\int_0^1 \frac{d\xi}{(\theta_0 + \theta_0^{-1}) + (\theta_0 - \theta_0^{-1})\cos(\pi\xi)}$$

is a rational number. In this case the overall gradient $\nabla_x U$ is also periodic with period $m\pi/k$ for some integer $m$. This completes the proof. □

Once we identified bounded solutions as periodic, the next step is to seek the period for each periodic orbit.

**Lemma 4.2.** *The period of each bounded orbit associated with (2.12)-(2.14) is given by*

(4.2) $$\overline{T} = \frac{2}{k}\int_{-\pi/2}^{\pi/2} \frac{d\theta}{\theta_0 + \theta_0^{-1} + (\theta_0^{-1} - \theta_0)\sin\theta}.$$

*Here $\theta_0 = \theta_0(\alpha) < 1$ is given by*

(4.3) $$\theta_0 := \frac{\sqrt{1 + 8kp_0} - 1}{\sqrt{1 + 8kp_0} + 1},$$

where $p_0$ is determined by the initial data

$$p_0(\alpha) = \frac{\sqrt{2\varphi_0 - \Gamma_0 - 4k^2}}{d_0^2 + (\sqrt{2\varphi_0 - 4k^2 - \Gamma_0} - 2k)^2}.$$

*Proof.* Due to the symmetry it suffices to compute the half period. The intersection points of the ellipse $V(\varphi, d) = V(\varphi_\pm, 0)$ with $d = 0$ can be written explicitly in terms of the initial data

(4.4) $$\varphi_- = \frac{2k}{\sqrt{C_0}}\theta_0, \quad \varphi_+ = \frac{2k}{\sqrt{C_0}}\theta_0^{-1}.$$

Using the equation

$$\frac{d\varphi}{dt} = -\varphi d,$$

where along the trajectory from $(\varphi_-, 0)$ to $(\varphi_+, 0)$ in the lower-half phase plane $(\varphi, d)$ we have

$$d = -\sqrt{\mathcal{D}_0\varphi - C_0\varphi^2 - 4k^2} = -\sqrt{C_0}\sqrt{(\varphi_+ - \varphi)(\varphi - \varphi_-)}.$$

Therefore

(4.5) $$\overline{T} = 2\int_{\varphi_-}^{\varphi_+} \frac{d\varphi}{-\varphi d} = \frac{2}{\sqrt{C_0}} \int_{\varphi_-}^{\varphi_+} \frac{ds}{s\sqrt{(\varphi_+ - s)(s - \varphi_-)}}.$$

Let $s = \frac{\varphi_- + \varphi_+}{2} + \frac{\varphi_+ - \varphi_-}{2}\tau$, and using the expression of $\varphi_-$ and $\varphi_+$ in (4.4) one has

$$\overline{T} = \frac{4}{\sqrt{C_0}} \int_{-1}^{1} \frac{d\tau}{[\varphi_- + \varphi_+ + (\varphi_+ - \varphi_-)\tau]\sqrt{1 - \tau^2}}$$

$$= \frac{4}{\sqrt{C_0}} \int_{-\pi/2}^{\pi/2} \frac{d\theta}{\varphi_- + \varphi_+ + (\varphi_+ - \varphi_-)sin\theta}$$

$$= \frac{2}{k} \int_{-\pi/2}^{\pi/2} \frac{d\theta}{\theta_0 + \theta_0^{-1} + (\theta_0^{-1} - \theta_0)sin\theta}.$$

This gives the desired result. $\square$

## 5. Flow Map

For the smooth flow we may further study the structure of the flow map. Assume $x = X_\alpha(t)$ is the flow map started at the initial position $\alpha$, then one has

$$\dot{X}_\alpha := \frac{dX_\alpha}{dt} = U(X_\alpha), \quad X_\alpha(0) = \alpha,$$

and momentum equation can be written as

$$\ddot{X}_\alpha = 2kJ\dot{X}_\alpha.$$

Integration once gives

$$\dot{X}_\alpha = U_0(\alpha) + 2kJ(X_\alpha - \alpha),$$

where $U_0(\alpha)$ is the initial velocity at location $\alpha$. The above equation leads to the flow map expression

$$X_\alpha(t) = \frac{1}{2k}J^{-1}e^{2kJt}U_0(\alpha) + \alpha - \frac{1}{2k}J^{-1}U_0(\alpha).$$

This flow map determines the unique smooth velocity field if and only if the indicator matrix,
$$\Gamma := \frac{\partial X_\alpha(t)}{\partial \alpha} = I - \frac{1}{2k}J^{-1}(I - e^{2kJt})\nabla_\alpha U_0,$$
remains nonsingular. Noting that $J^{-1} = -J$ and
$$e^{2kJt} = \begin{pmatrix} cos(2kt) & sin(2kt) \\ -sin(2kt) & cos(2kt) \end{pmatrix},$$
we have
$$\Gamma = I + \frac{1}{2k}J\begin{pmatrix} 1 - cos(2kt) & -sin(2kt) \\ sin(2kt) & 1 - cos(2kt) \end{pmatrix}\nabla_\alpha U_0.$$
Hence, with $U = (u, v)$ we find
$$2k\Gamma = 2kI + \begin{pmatrix} sin(2kt) & 1 - cos(2kt) \\ cos(2kt) - 1 & sin(2kt) \end{pmatrix}\nabla_\alpha U_0$$
$$= \begin{pmatrix} 2k + v_{0x} + u_x sin(2kt) - v_{0x}cos(2kt) & v_{0y} + u_{0y}sin(2kt) - v_{0y}cos(2kt) \\ -u_{0x} + v_{0x}sin(2kt) + u_{0x}cos(2kt) & k - u_{0y} + v_{0y}sin(2kt) + u_{0y}cos(2kt) \end{pmatrix}.$$

A careful calculation gives its determinant as
$$det(2k\Gamma) = 4k^2 - 2k\omega_0 + 2det(\nabla_\alpha U_0) + (2k\omega_0 - 2det(\nabla_\alpha U_0))cos(2kt) + (2kd)sin(2kt).$$
Thus $\Gamma(t)$ remains nonsingular for all time if and only if $det(2k\Gamma) \neq 0$, i.e.,
$$4k^2 - 2k\omega_0 + 2det(\nabla_\alpha U_0)$$
$$\notin \left(-\sqrt{(2k\omega_0 - 2det(\nabla U_0))^2 + 4k^2 d_0^2}, \sqrt{(2k\omega_0 - 2det(\nabla_\alpha U_0))^2 + 4k^2 d_0^2}\right),$$
which is equivalent to
$$(5.1) \qquad (4k^2 - 2k\omega_0 + 2det(\nabla_\alpha U_0))^2 > (2k\omega_0 - 2det(\nabla_\alpha U_0))^2 + 4k^2 d_0^2.$$
The spectral gap $\Gamma_0$ relates the determinant and the divergence as
$$\Gamma_0 = d_0^2 - 4det(\nabla_\alpha U_0).$$
This when applied to the above inequality (5.1) gives
$$\Gamma_0 < 4k^2 - 4k\omega_0,$$
which is exactly the critical threshold (1.5) in phase space stated in Theorem 1.1.

## 6. Kinetic formulation

This section describes a kinetic formulation for the rotational model (1.2) in terms of a density function, $f = f(t, x, \xi)$ governed by the BGK model
$$(6.1) \qquad \partial_t f + \xi \cdot \nabla_x f + 2k\xi^\perp \cdot \nabla_\xi f = \frac{1}{\epsilon}(M - f), \quad \xi^\perp := J\xi.$$
Here $M = M_{\{\rho, U\}}(\xi)$ is the Maxwellian given by
$$M = \frac{\rho}{\sqrt{\pi T}}e^{-|\xi - U|^2/T}, \quad \xi = (\xi_1, \xi_2) \in \mathrm{I\!R}^2.$$

The fixed temperature $T$, plays no role in this pressureless model. The connection between the distribution function $f$ and macroscopic flow variable is realized in terms of the usual moments of density $\rho$, momentum $m = \rho U$ and total energy $E = \rho |U|^2/2$,

$$(\rho, \rho U, E)^\top = \int \psi^\top(\xi) f d\xi, \quad \psi(\xi) := \left(1, \xi, \frac{|\xi|^2}{2}\right)^\top.$$

The conservation principle for mass, momentum and energy during the course of particle collisions requires the equilibrium to satisfy the compatibility condition

$$\int (M - f)\psi^\top(\xi) d\xi = 0,$$

while the rotational forcing is introduced through the potential

$$\int 2k\xi^\perp \cdot \nabla_\xi f \psi^\top(\xi) d\xi = (0, -2k\rho JU, 0).$$

Indeed, a straightforward computation yields $\int 2k\xi^\top \cdot \nabla_\xi f d\xi = 0$; for the momentum equation we compute

$$\int 2k\xi^\perp \cdot \nabla_\xi f \xi d\xi = 2k \int \xi(\xi_2 \partial_{\xi_1} f - \xi_1 \partial_{\xi_2} f) d\xi$$
$$= 2k \int \begin{pmatrix} \xi_1 \xi_2 \partial_{\xi_1} f \\ -\xi_1 \xi_2 \partial_{\xi_2} f \end{pmatrix} d\xi$$
$$= 2k \begin{pmatrix} -\rho v \\ \rho u \end{pmatrix}$$
$$= -2k\rho JU,$$

while the presence of this forcing does not change the energy equation since

$$2k \int \xi^\perp \cdot \nabla_\xi f |\xi|^2 d\xi = 0.$$

The first three moments of (6.1) then yield the equivalent extended system of (1.3),(1.4),

$$\partial_t \begin{pmatrix} \rho \\ \rho U \\ \rho \frac{|U|^2}{2} \end{pmatrix} + \nabla_x \cdot \begin{pmatrix} F_\rho \\ F_m \\ F_E \end{pmatrix} = \begin{pmatrix} 0 \\ 2k\rho JU \\ 0 \end{pmatrix}.$$

The corresponding macroscopic fluxes are

$$(F_\rho, F_m, F_E)^\top := \int \psi^\top(\xi) \xi f d\xi,$$

and under the closure $f = M_{\{\rho, U\}}$ we conclude

$$(F_\rho, F_m, F_E)^\top = \left(\rho U, \rho U \otimes U, \rho \frac{|U|^2}{2}\right)^\top.$$

**ACKNOWLEDGMENTS:** Research was supported in part by ONR Grant No. N00014-91-J-1076 (ET) and by NSF grant #DMS01-07917(ET, HL).

Department of Mathematics, Iowa State University, Ames, IA 50011.
*E-mail address*: `hliu@iastate.edu`

Department of Mathematics, Center of Scientific Computation And Mathematical Modeling (CSCAMM) and Institute for Physical Science and Technology (IPST), University of Maryland, MD 20742.
*E-mail address*: `tadmor@cscamm.umd.edu`